\documentclass[12pt,a4paper]{amsart}
\usepackage[margin=3cm]{geometry}
\title{Local automorphisms of the Hilbert ball}
\author{Bernhard Lamel}
\address{Universit\"at Wien, Fakult\"at f\"ur Mathematik \\ Nordbergstrasse 15, A-1090 Wien, \"Osterreich}
\email{lamelb@member.ams.org}%
\thanks{The author was supported by the FWF, Projekt P17111}

\newcommand{\inp}[2]{\langle #1 , #2 \rangle}

\newcommand{\vnorm}[1]{\left\|  #1 \right\|}

\newcommand{\C}{\mathbb{C}}
\newcommand{\CN}{\mathbb{C}^N}

\newcommand{\R}{\mathbb{R}}

\newcommand{\N}{\mathbb{N}}

\newcommand{\dop}[1]{\frac{\partial}{\partial #1}}

\newcommand{\hs}{\mathsf{H}}
\newcommand{\hb}{\mathsf{B}}
\newcommand{\ball}{\mathbb{B}}

\newlength{\extendaxesby}

\DeclareMathOperator{\imag}{Im}

\newtheorem{thm}{Theorem}

\theoremstyle{definition}

\newtheorem{exa}{Example}
\newtheorem{rem}{Remark}

\begin{document}
\begin{abstract}{We prove an analogue of Alexander's Theorem \cite{Al1} for holomorphic mappings
of the unit ball in a complex Hilbert space: Every holomorphic mapping which 
takes a piece of the boundary of the unit ball into the boundary of the 
unit ball and whose differential at
some point of this boundary is onto is the restriction of an automorphism of the 
ball. We also show that it is enough to assume that the mapping is only G\^ateaux-holomorphic.}
\end{abstract}
\maketitle
\section{Introduction and Statement of Results}
Let $\ball_N\subset\CN$ be the unit ball in $\CN$, that
is, 
\[ \ball_N = \left\{ Z \in \CN \colon 
\vnorm{Z} < 1\right\}.\]
It is well known that if $H\colon(U,p_0) \to \CN$ is a 
holmorphic mapping
defined in a neighbourhood $U\subset\CN$ 
of some $p_0\in\partial\ball_N$
which maps $U\cap \partial\ball_N$ into $
\partial\ball_N$, 
then $H$ is either 
a constant mapping, or it extends to an automorphism of $\ball_N$.
This Theorem is due to Alexander \cite{Al1}. 

In a Hilbert space of infinite dimension, Alexander's Theorem is not 
valid (see Example~\ref{e:homogsum}; the author is 
grateful to Laszlo Lempert for pointing out this fact). However, if we 
assume in addition that the differential of the 
mapping $H$ (see
Theorem~\ref{t:mainball}) is {\em onto} at 
some point $p_0$ on the boundeary of the ball, then 
the same extension property as in finite dimensional Hilbert space holds. 
That is, a mapping $H$
defined near a point $p_0$ in  the boundary of the
unit ball $\hb \subset \hs$
of some complex Hilbert space $\hs$ whose differential at $p$
is onto extends to an automorphism of $\hb$. In order to formulate our
result, we recall that a mapping $H\colon U \to \hs$ 
defined on some open subset
$U\subset \hs$ is {\em G\^ateaux-holomorphic} if the function
$Z\mapsto \inp{H(Z)}{\tilde Z}$ is 
holomorphic on all sets of the form $E\cap U$ where 
$E$ is some finite dimensional, affine subspace of $\hs$. A
G\^ateaux-holomorphic 
map $H$ is {\em holomorphic } if it is also continuous (for these definitions,
see e.g. the book of Dineen \cite{Dineenbook}).
We can now state our main result.

\begin{thm}
  \label{t:mainball} Let $\hs$ be a complex Hilbert space, 
  $P\in\hs$ with $\vnorm{P}=1$, and
  $H \colon U \to \hs$ be a G\^ateaux-holomorphic 
  map defined on a neigbourhood 
  $U$ of $P$ with the property that $\vnorm{H(Z)} = 1$ for 
  all $Z\in U$ with $\vnorm{Z} =1$, and such 
  that the range of $H^{\prime} (P)$ is $\hs$. 
  Then $H$ is holomorphic, and furthermore, 
  $H$ extends to an automorphism of $\hb$.
\end{thm}

Note that we do not assume any continuity of the map here; it 
is indeed part of the conclusion that the mapping $H$ is
continuous. Also note the assumption of G\^ateaux-holomorphicity
implies that $H$ has 
a weak (not {\em a priori} continuous) derivative, which
is defined on all of $\hs$; it is this 
derivative that we denote by $H'(P)$ in the statement 
of Theorem~\ref{t:mainball}. It follows that 
$H'(P)v$ is given by the strong limit of the difference quotients, 
\[ H'(P) v = \lim_{\C\ni t\to 0} \frac{H(P+tv) - H(P)}{t}, \]
we do not assume {\em a priori} that $H'(P)$ is bounded; 
similar remarks hold for higher derivatives. 
We note that the automorphisms of a ball in 
Hilbert space are well known (this goes back 
to a paper of Renaud \cite{Re1}), and that any such automorphism 
is a linear fractional map (just as in the finite dimensional case) which
extends holomorphically across the boundary of the ball. 

The proof of Theorem~\ref{t:mainball} is an application of 
the technique of Segre varieties 
on the Hilbert space in question. The special
structure of the boundary of the ball allows us 
to do the necessary polarization very explicitly 
and also 
allows us to exploit the Riesz Representation Theorem
(that is, the self-duality of the Hilbert space)
to explicitly calculate expressions for $H$. We will do 
this computation in a set of different coordinates, which 
are used for computations of 
this type in finite dimensional 
complex space, the so-called ,,normal 
coordinates''; in the case of the 
unit ball, its boundary is put into normal coordinates by 
the well-known Cayley transform. 
Our calculations are very explicit, and
also yield a good understanding of the 
structure of the automorphism group in question.

In our proof, we will pass back and forth between 
applying the Segre technique along finite dimensional
affine subspaces and using the results of these
calculations to get global expressions for the map. 
It is this passing back and 
forth that allows us to establish holomorphicity of the
map without assuming continuity. The technique of the 
Segre varieties has been used extensively in the study of 
holomorphic mappings in the finite-dimensional case; for 
a detailed discussion of their use in that case, 
we refer the reader to e.g. the 
survey article of Baouendi, Ebenfelt and Rothschild \cite{BERbull}.

Algorithmically, our calculations follow the path along which
Baouendi, Ebenfelt and Rothschild \cite{BER3} established that automorphisms
of real-analytic hypersurfaces in $\CN$ are parametrized by their
jets of some finite order. Our calculations here also yield 
an explicit parametrization for the group  of local automorphisms of 
the boundary of the unit ball in $\hs$; see \eqref{e:par} below.

As noted above, the change to normal coordinates is given by means of the 
Cayley transform, which is recalled in Section~\ref{s:coord}.
In these coordinates, given a point $P\in\partial \hb$, 
$\hb$ corresponds to 
\[ \mathbb{H}_{+} = \left\{(z,w)\colon \imag w > \vnorm{z}^{2} \right\} \subset
\C \times P^{\perp},  \]
and $P$ corresponds to $(0,0)$. The boundary $\partial \hb$ 
corresponds to 
\[ \mathbb{H} = \left\{ (z,w)\colon \imag w = \vnorm{z}^2 \right\}. \]
Theorem~\ref{t:mainball} is a consequence of the following, more 
detailed Theorem: 
\begin{thm}\label{t:mainh}
  Let $F$ be a Hilbert space, $H\colon U\to
  F\times \C$ a map defined in a neighbourhood
  $U$ of $(0,0)$ in $F\times\C$ which is 
  G\^ateaux-holomorphic on $U$ (that is, for every
  affine subspace $E$ of $F\times\C$, and for
  every continuous linear form $\phi $ on $F\times \C$,
  the map 
  $\phi (H) \colon E\cap U \to \C$ is holomorphic), 
  satisfies $H(0) = 0$,
  and $H(U\cap \mathbb{H})\subset \mathbb{H}$. Assume in addition
  that  the range of $H^\prime (0)$ is $\hs$. Then 
  $H$ is holomorphic, and furthermore, if we
  denote the coordinate in $F$ by $z$ and the 
  coordinate in $\C$ by $w$, then 
  \begin{equation}
    H(Z) = \left( 
      s U \frac{ z + a w }{1- 2i \inp{z}{a} 
      + (R - i \vnorm{a}^2) w}
      ,s^2  \frac{  w }{1- 2i \inp{z}{a} 
      + (R - i \vnorm{a}^2) w}
    \right)
    \label{e:Hform}
  \end{equation}
  where $U$ is a unitary map, $s\in\R_{+}$, $a\in \C$,
  $R\in\R$, and $U$, $s$, $a$, $R$ are given by
  \begin{equation} 
    \begin{aligned}
      U &= \frac{f_z (0)}{\sqrt{g_w (0)}}, &
      s &= {\sqrt{g_w (0)}}, \\
      a &= -  f_z (0)^{-1} f_w (0) , &
      R &= \frac{2 g_{w^2 } (0) + i \vnorm{f_{w}(0)}^2}{g_{w} (0)}.
    \end{aligned}\label{e:par}
  \end{equation}
\end{thm} 

It is natural to ask whether or not the assumption that the range
of $H^\prime (P)$ is necessary. The simplest example of a map violating
this condition, which is not an automorphism of $\hb$ was pointed out 
to the author by Lempert: Just consider an isometric map of $\hs$ onto
a proper subspace of itself (e.g. the shift map on $\ell^2$).  

There are numerous other examples of maps whose derivative is not onto 
at any point on the boundary which are interesting; let us 
discuss two of them here.

\begin{exa}\label{e:homogsum}
  Let $\hs = \ell^2$, and write $Z = (z_1,z_2,\cdots)$ for the variable in $\ell^2$.
  Consider the set 
  \[ \Gamma = \bigcup_{k\in\N} \N^k.\]
  This is a countable set; choose a bijection $\phi:\Gamma\to\N$. Let
  $\lambda_j$ be a sequence with 
  \[ \sum_{j} |\lambda_j|^2 = 1.\]
  For $\alpha = (\alpha_1,\cdots,\alpha_k) \in \Gamma$, write
  $z_\alpha = z_{\alpha_1} \cdots z_{\alpha_n}$ and $|\alpha| = k$.
  We define $H(Z) = W = (w_1,w_2,\cdots)$ by setting $w_{\phi(\alpha)}
  = \lambda_{|\alpha|} z_{\alpha}$. Then $H$ maps $r \hb$ into
  $R(r) \hb$ for $0\leq r < r_0$ for some $r_0 > 1$, where
  \[ R(r) = \sum_j |\lambda|^2 r^j. \]
  It is easy to show that $H$ is a holomorphic map; it is clearly not
  an automorphism. 
\end{exa}

\begin{exa}
  Our second example generalizes the Whitney map. Using the notation introduced
  above, for a given $p \in \N\cup\left\{ \infty \right\}$,
  we now consider as an index set
  \[ \Gamma_p = \bigcup_{ 0\leq q < p} \left\{ q \right\} \times \left\{ 
  j\in\N \colon j > 1\right\} \cup \{p\}, \quad p< \infty,\]
  \[ \Gamma_\infty = \bigcup_{ 0\leq q } \left\{ q \right\} \times \left\{ 
  j\in\N \colon j > 1\right\} .\]
  Again, we choose a bijection $\phi \colon \Gamma_p \to \N$ and 
  we define $W = H(Z)$ by 
  \[ w_{\phi(q,k)} = z_1^q z_k, \quad w_{\phi(p)} = z_1^p 
  \text{ in case } p<\infty.\] 
  Then we have that 
  \[ \vnorm{H(Z)}^2= \frac{1- |z_1|^{2p}}{1 - |z_1|^2} 
  \left( \vnorm{Z}^2 - |z_1|^2 \right) + |z_1|^{2p}, \quad p \neq \infty,\]
  \[ \vnorm{H(Z)}^2 = \frac{1}{1-|z_1|^2} \left( \vnorm{Z}^2
  - |z_1|^2\right), \quad p = \infty.\]
  It is easy to check that if $p\neq \infty$, $H$ extends holomorphically to a 
  neighbourhood of $\overline{\hb}$. On the other hand, if $p=\infty$, 
  $H$ {\em does not extend} holomorphically to any neighbourhood of $(1,0,\cdots)$
  (indeed, it does even extend continuously up to that point!). 
\end{exa}
\section{Changing coordinates and some observations}\label{s:coord}

We now assume that we have the following situation: 
$H$ is a map defined in a neighbourhood $U$ of $P$, satisfying 
the assumptions of Theorem~\ref{t:mainball}. After composing 
with a rotation in the plane spanned by $P$ and $H(P)$, we 
can assume that $H(P) = P$. We decompose $\hs = F \oplus 
\C P$, where $F = P^\perp$. In this decomposition, 
we write $Z = (\zeta, \eta)$ with $\zeta \in F$ and
$\eta \in \C$. Our goal in this section is to show that
Theorem~\ref{t:mainh} implies Theorem~\ref{t:mainball},
and to assemble the necessary prerequisites for the proof of 
Theorem~\ref{t:mainh} which will be given in the following
section.

\subsection{The Cayley transform}
The Cayley transform (resp. its inverse) is defined by 
\[
  (z,w) = \left(\frac{\zeta}{1+\eta}, 
  i \frac{1-\eta}{1+\eta}\right), \quad  
  ( \zeta,\eta)= \left( \frac{2iz}{i+w},
 \frac{i-w}{i+w}\right),
\]
and is a local biholomorphism from a neigbourhood of $P$ to
$F\oplus\C$, taking $P$ to $(0,0)$ and 
the boundary of the unit ball (with the point $-P$ omitted) to 
the hypersurface $\mathbb{H}$ in $F\oplus \C$ defined by the equation
\[ \imag w = \vnorm{z}^2, \quad (z,w) \in F\oplus \C.\]
It also 
takes the interior of the unit ball to the
half-space $\mathbb{H}_{+} = \{ \imag w \geq \vnorm{z}^2 \}$. 

Employing the Cayley transform, we see that Theorem~\ref{t:mainh} implies
Theorem~\ref{t:mainball}. So from now on, we assume that
$H$ is defined in a neighbourhood of $(0,0)\in F\oplus\C$,
takes $(0,0)$ to $(0,0)$, and
if we write $H(z,w) = (f(z,w),g(z,w))$ with $f(z,w) \in F$ and
$g(z,w) \in \C$, then we have
\[ \imag g(z,w) = \vnorm{f(z,w)}^2, \text{ if } 
\imag w = \vnorm{z}^2.\]

\subsection{Complexification}
We will now complexify the last equation 
(in a bit of a nonstandard
manner, owing to the fact that 
the natural isomorphism $\hs \to \hs^*$
induced by 
the inner product is conjugate linear). We claim that we have
\begin{equation}
  g(z,w) - \overline{g(\chi,\tau)} = 2i
  \inp{f(z,w)}{f(\chi,\tau)}, \text{ if } w - \bar \tau = 
  2i \inp{z}{\chi}.
  \label{e:b1}
\end{equation}
We check this claim by restricting to the 
finite dimensional subspace spanned by $z$ and $\chi$: 
Choose an orthonormal basis $\{ u, v \}$ of the span of 
$z$ and $\chi$; We have
\[ z = r_z u + s_z v, \quad \chi = r_\chi u
+ s_\chi v, \]
and consider the hypersurface $M$ in $\C^3$ with 
coordinates $(r,s,w)$ defined by
\[ \imag w = |r|^2 + |s|^2.\]
On this hypersurface, consider the function $\rho$ 
defined  by 
\[ \rho (r,s,w) = \imag g(r u + s v,w) 
- \inp{f(r u + s v,w)}{f(z,w)}, \quad z = r u + s v. \]
A fact which we are going to use often in
the following is that the Banach-Steinhaus 
Theorem guarantees that a mapping which is 
G\^ateaux-holomorphic is automatically uniformely bounded 
along compact subsets of finitely dimensional affine subspaces
of $\hs$. This in turn implies that 
the function $\rho $ is real-analytic on a neighbourhood
of $0$ in $\C^3$; by assumption it vanishes on $M$; so
the usual complexification
shows that
\[ g(r u + s v,w) - 
\overline{g(p u + q v,\tau)} - 
2i \inp{f(r u + s v,w)}{f(p u + q v,\tau)} 
\text{ if } w - \bar \tau = 
2i ( r \bar p + s \bar q ) . \]

\subsection{CR and transversal
vector fields} 
Our next step is to introduce two vector fields which we
are going to use to differentiate \eqref{e:b1}. We choose an 
$u \in F$, and define 
\begin{equation}
  \mathcal{L}_u = \dop{ \bar u} - 2i \inp{z}{u} \dop{\bar\tau}.
  \label{e:Ldef}
\end{equation}
Application
of $\mathcal{L}_u$  to a function $\phi$ 
is defined in the following natural way: Suppose that 
$\phi (z,w,\chi,\tau)$ is valued
in some Hilbert space $E$, defined in a neighbourhood of 
$0$ and that we are
given $z_0$, $w_0$, $\chi_0$ and $\tau_0$. 
We choose an orthonormal basis $\{u,g,h\}$ of 
the space spanned by $z_0$, $\chi_0$ and $u$. Setting 
$ z = r u + s g + t h$, and $\chi = a u +  b g + ch$, we 
define $\mathcal{L}_u \phi (z_0,w_0,\chi_0,\tau_0)$ to
be the application of the vector field
\[ \dop{\bar a} - 2i r \dop{\bar \tau} \]
to $\phi$ at the appropriate point.

We also define the vector field 
\[ \mathcal{X} = \dop{w} + \dop{\bar \tau}.\]

\begin{rem}\label{rem:weakdiff}
  Note that the Banach-Steinhaus Theorem implies that 
given a map $h$ which is a G\^ateaux-holomorphic map
from a Hilbert space $\hs_1$ into some other Hilbert space
$\hs_2$, 
$h$ possesses all weak directional derivatives.
By this we mean the following: for each $u \in \hs_1$ and
each $u_0$ in the domain of $h$, there exists a vector
$v\in\hs_2$ such that for all $P \in \hs_2$, 
\[ \lim_{t\to 0} \frac{\inp{h(u_0 + t u)}{P} - \inp{h(u_0)}{P}}{t}
= \inp{v}{P}.\]
It follows by well-known arguments that this limit is actually strong, that
is, 
\[ v= \lim_{t\to 0} \frac{h(u_0 + tu) - h(u_0)}{t}.\]
In this situation, we will write $v= h'(u_0)(u)$. Note that
we do not assume that $h$ is differentiable when we use this 
notation; also, $h'(u_0)$ is not necessarily bounded.
In the case that we are dealing with, we will write
\[ f_z (z_0,w_0) u =
f'(z_0,w_0)(u,0), \]
and so on; we will identify e.g. $f_w (z_0,w_0)$ with the
vector $f_w(z_0,w_0)(0,1)$ and the derivatives $g_w(z_0,w_0)$ 
with number $g_w (z_0,w_0)(0,1)$. Let us emphasize again
that writing $f_z$ does {\em not} mean we assume 
that $f$ is differentiable (that's part of what we want to show).
\end{rem}

\section{Proof of Theorem~\ref{t:mainh}}
We now start by applying $\mathcal{L}_u$ to \eqref{e:b1}. We
obtain
\begin{equation}
  \begin{aligned}
    - \overline{ g_z (\chi,\tau) u + 
    2 i \inp{u}{z} g_w (\chi,\tau) } &= 2i 
    \inp{f(z,w)}{f_z (\chi,\tau) u + 2 i \inp{u}{z} 
    f_\tau (\chi,\tau)},
  \end{aligned}\label{e:Ltob1}
\end{equation}
which is valid for $w - \bar \tau = 2i \inp{z}{\chi}$ and 
every $u \in F$. Setting $\chi = 0$, $ \tau = w = 0$ and noting that
\eqref{e:b1} implies $g(z,0) = 0$, we
get that
\begin{equation}
  \overline{g_w(0)} \inp{z}{u} = \inp{ f(z,0) }{ f_z (0) u + 
  2i \inp{u}{z} f_w (0)}
  \label{e:levi1}
\end{equation}
which holds for every $z\in F$ close 
enough by $0$ and every $u\in F$. In particular, we have
\begin{equation}
  \overline{g_w (0)} \inp{u}{v} = \inp{f_z (0) u}{f_z(0) v}.
  \label{e:levi2}
\end{equation}
So, if the range of $H^\prime (P)$ is $\hs$, this 
implies that $f_z (0)$ is an isomorphism. Furthermore, in that
case,
$g_w (0) = r^2$ is a positive real number, and the map 
$r^{-1} f_z (0)$ is unitary. Consider the
automorphisms $\omega_{U,s}$ for unitary
operators $U$ of $F$ and $s \in \R$ 
defined by 
\[ \omega_{U,s} (z,w) = (sUz,s^2 w).\]
Composing with 
the automorphism 
$(z,w) \mapsto \omega_{r f_z(0)^{-1}, r^{-1}} (z,w) = 
(f_z (0)^{-1}  z , r^{-2} w)$,
we can assume that $f_z (0) = I$ and $g_w (0) = 1$. 

Returning to \eqref{e:levi1}, we note that we can now 
explicitly compute $f(z,0)$. Indeed, with the 
simplifications made above, we
have
\begin{equation}
  \inp{z}{u} = \inp{f(z,0)}{u + 2 i \inp{u}{z} f_w(0)}.
  \label{e:firstsegre1}
\end{equation}
vor $z\in F$ close by $0$, and all $u\in F$.
Writing $T u  = u + 2i \inp{u}{z} f_w (0)$, we
conclude that $f (z,0) = T^{* -1} z$, and computing 
$T^{*-1} v$ we obtain
\[ f(z,0) = \frac{z}{1 + \inp{z}{2i f_w (0)}}.\]

The next simplification we can make is that we can 
assume $f_w (0) = 0$. 
Indeed, for $a \in F$, the automorphism 
$\phi_a$ defined by 
\[\phi_a (z,w) = 
\left( \frac{z + wa}{1 - 2i \inp{z}{a} - i w \vnorm{a}^2},
\frac{w}{1 - 2i \inp{z}{a} - i w \vnorm{a}^2}\right)\]
takes $\mathbb{H}$ into itself. 
Composing with $\phi_a$ for $a = -f_w (0)$ we see that
we can actually assume $f_w (0)=0$. 

In the next step, we first apply $\mathcal{X}$ to \eqref{e:b1}
to obtain
\begin{equation} g_w (z,w) - \overline{g_w (\chi,\tau)} =
  2i \left( \inp{f_w(z,w)}{f (\chi,\tau)} + 
  \inp{f(z,w)}{f_w (\chi,\tau)} \right)\label{e:xtob1}
\end{equation}
Evaluating again at $\chi=0$, $\tau = w =0$,
we get $g_w (z,0) = \overline{g_w (0)}$.
Now we apply $\mathcal{L}_u$ to \eqref{e:xtob1},
and evaluate at $\chi = 0$, $w = \tau =0$. 
The result is 
\begin{equation}
  \overline{g_{w^2} (0)} \inp{z}{u} = \inp{f_w(z,0)}{u } 
  + \inp{z}{f_{zw} (0) u + 2i \inp{u}{z} f_{w^2} (0)},
  \label{e:fwseg1}
\end{equation}
from which we see that $f_{zw} (0)^{*} z$ is defined
for all $z\in F$ (since the equation
exhibits that $u\mapsto \inp{z}{f_{z,w} (0) u}$ is continuous for
$z$ close by $0$).
We have
\begin{equation}
  f_{w}(z,0) = \overline{g_{w^2} (0)} z
  - f_{zw}^{*} (0) z 
  - \inp{z}{2if_{w^2} (0)} z
  \label{e:fwiseq}
\end{equation}
Now turn back to \eqref{e:Ltob1} and evaluate this equation at
$\tau = 0$, $w = 2i \inp{z}{\chi}$. Applying
our simplifications, it then reads
\begin{equation}\begin{aligned}
  \inp{z}{u} &
  =\overline{g_w (0)} \inp{z}{u} \\
  &=\overline{g_w (\chi,0)} \inp{z}{u}\\ 
  &= \inp{f(z,2i \inp{z}{\chi})}{f_z (\chi,0)u + 2i \inp{u}{z}
  f_{w}(\chi,0)} \\
  &= \inp{f(z,2i \inp{z}{\chi})}{ u + 2i \inp{u}{z}
  f_{w} (\chi,0) }.
  \label{e:secseg}
\end{aligned}
\end{equation}
We note that $Tu = u+ 2i \inp{u}{z} f_{w} (\chi,0)$ 
is a continuous, invertible linear map for $z$ and $\chi$
close by $0$, with and compute 
\[ T^{*-1} v = v - \frac{ \inp{v}{2i f_w (\chi,0)}}
{1 + \inp{z}{2i f_w (\chi,0)} }z;\]
with this observation, \eqref{e:secseg} implies that
\begin{equation*}\begin{aligned}
  f(z,2i \inp{z}{\chi}) &= 
  T^{*-1} z \\ 
  &= \frac{z}{1 + \inp{z}{2i f_{w } (\chi,0)}}\\
  &=
  \frac{z}{ 1 - 2i g_{w^2} (0) \inp{z}{\chi} 
  + 2i \inp{z}{f_{zw}(0)^{*} \chi }
  - 4 \inp{z}{\chi} \inp{f_{w^2}(0)}{\chi}},
\end{aligned}
\end{equation*}
and
\[ g(z,2i \inp{z}{\chi} ) = 
\frac{2i\inp{z}{\chi}}{ 1 - 2i g_{w^2} (0) \inp{z}{\chi} 
+ 2i \inp{z}{f_{zw}(0)^{*}  \chi }
- 4 \inp{z}{\chi} \inp{f_{w^2}(0)}{\chi}}.\]
We now claim that $f_{w^2}(0) = 0$. To check this, 
we have to use G\^ateaux-holomorphicity once more.
We choose a $\chi_0\in F$ and
for $z$ with $\inp{z}{\chi_0} \neq 0$
and we set $w = 2i \inp{z}{\lambda \chi_0}$, getting 
$\lambda =  \frac{i\bar w}{2 \inp{\chi_0}{z} } $.
We also set 
$\inp{z}{f_{zw} (0)^{*} \chi_0} = a(\chi_0) \inp{z}{\chi_0}$.
\begin{equation}\begin{aligned}
  f(z,w) &= \frac{z}{1 - (g_{w^2} (0) -  a(\chi_0)) w -
  w^{2} \frac{\inp{f_{w^2}(0)}{\chi_{0}}}{\inp{z}{\chi_0}}}\\
  g(z,w) &= \frac{w}{1 - (g_{w^2} (0) -  a(\chi_0)) w -
  w^{2} \frac{\inp{f_{w^2}(0)}{\chi_{0}}}{\inp{z}{\chi_0}}}\\
  \label{e:wpar}
\end{aligned}
\end{equation}
Since $H$ is G\^ateaux-holomorphic, the map 
\[ \psi (\lambda,\mu) = g(\lambda \chi_0, \mu) = 
\frac{\mu}{ {1 - (g_{w^2} (0) -  a(\chi_0)) \mu -
\frac{\mu^{2}}{\lambda} \frac{\inp{f_{w^2}(0)}{\chi_{0}}}{\vnorm{\chi_0}^2}}}\]
is a holomorphic map $\C^2 \to \C$ (a priori, 
the equation only holds on an open subset of $\C^2$, of 
course, but this implies that the functions coincide
as rational functions). This of course is 
only possible if the coefficient of $\frac{\mu^2}{\lambda}$
in the denominator does vanish, that is, if $\inp{f_{w^2} (0)}
{\chi_0} = 0$. Since this holds for all $\chi_0$ sufficiently 
close by $0$, we conclude that $f_{w^2} (0) = 0$.
Thus, \eqref{e:wpar} becomes 
\begin{equation}\begin{aligned}
  f(z,w) &= \frac{z}{1 - (g_{w^2} (0) -  a(\chi_0)) w }\\
  g(z,w) &= \frac{w}{1 - (g_{w^2} (0) -  a(\chi_0)) w }
  \label{e:wpar2}
\end{aligned}
\end{equation}
This holds
for all $(z,w)$ close 
by the origin which are of the form $w = 2i \inp{z}{\chi}$
for some $\chi$ close by the origin. Since 
this definitely includes an open set in each
finite dimensional affine subspace, the 
principle of analytic continuation and G\^ateaux-holomorphicity
of $H$ imply that \eqref{e:wpar2} holds
 for all $(z,w)$ close by the origin.
At this point, 
we have thus established that $H$ is actually holomorphic
in a neighbourhood of the origin. 
Taking the second derivative with respect to $w$ of 
the equation for $g$ in \eqref{e:wpar2}, we see that 
$ a(\chi_0) = \frac{g_{w^2} (0)}{2} $, so that
$ f_{zw} (0) u = \frac{g_{w^{2}} (0)}{2} u$, 
and we conclude that 
\begin{equation}\begin{aligned}
  f(z,w) &= \frac{z}{1 + R  w }\\
  g(z,w) &= \frac{w}{1  + R w  }
  \label{e:wpar3}
\end{aligned}
\end{equation}
where $R = - \frac{g_{w^{2}} (0)}{2}$. Next,
let us check 
that $H$ leaving $\mathbb{H}$ invariant
implies that $R\in \R$. We have
\begin{equation}
  \begin{aligned}
    \imag \frac{w}{1+ Rw} &= \frac{ \imag w }{|1+ Rw|^2} 
    - \frac{ |w|^2 (\imag R)}{|1+ Rw|^2}\\
    &= \frac{\vnorm{z}^2}{|1+Rw|^2}
  \end{aligned}
  \label{e:checkgw2inR}
\end{equation}
which implies that $\imag R = 0$. So, we have
$g_{w^2} (0) \in \R$. Thus, the 
denominator of our expressions for $f$ and $g$ does
not vanish on the set $ \mathbb{H}_{+} =
\{ \imag w \geq {\vnorm{z}^2} \} $; 
thus
an  $H$ of the form above 
is holomorphic on all of $\mathbb{H}_{+}$.

We will now recover the general form of a map $H$ by retracing
the simplifications made above. Let us write $H_R$ for the
map given by \eqref{e:wpar3}. Then we have to compose with 
a map of the form $\phi_a$ and a map of the form
$\omega_{U,s}$. 
This yields as the general form of an automorphism
\begin{equation}
  \omega_{U,s} \phi_a H_R = \left( 
  s U \frac{ z + a w }{1- 2i \inp{z}{a} + (R - i \vnorm{a}^2) w}
  ,s^2  \frac{  w }{1- 2i \inp{z}{a} + (R - i \vnorm{a}^2) w}
  \right)
\end{equation}
We are only missing the last conclusion of Theorem~\ref{t:mainh}; in 
order to prove this, just note that the particular
choice of the maps made above implies that
we can compute $U$, $s$, $ a$ and $R$ from the
derivatives of 
our original map $H$ in the following way:
\[ \begin{aligned}
  U &= \frac{f_z (0)}{\sqrt{g_w (0)}}, &
  s &= {\sqrt{g_w (0)}}, \\
  a &= -  f_z (0)^{-1} f_w (0) , &
  R &= \frac{2 g_{w^2 } (0) + i \vnorm{f_{w}(0)}^2}{g_{w} (0)}.
\end{aligned}\]

\bibliographystyle{plain}
\bibliography{bibliography}

\begin{thebibliography}{1}

\bibitem{Al1}
H.~Alexander.
\newblock Holomorphic mappings from the ball and polydisc.
\newblock {\em Math. Ann.}, 209:249--256, 1974.

\bibitem{BER3}
M.~S. Baouendi, P.~Ebenfelt, and L.~P. Rothschild.
\newblock Parametrization of local biholomorphisms of real analytic
  hypersurfaces.
\newblock {\em Asian J. Math.}, 1(1):1--16, 1997.

\bibitem{BERbull}
M.~S. Baouendi, P.~Ebenfelt, and L.~P. Rothschild.
\newblock Local geometric properties of real submanifolds in complex space.
\newblock {\em Bull. Amer. Math. Soc. (N.S.)}, 37(3):309--336 (electronic),
  2000.

\bibitem{Dineenbook}
S.~Dineen.
\newblock {\em Complex analysis on infinite-dimensional spaces}.
\newblock Springer Monographs in Mathematics. Springer-Verlag London Ltd.,
  London, 1999.

\bibitem{Re1}
A.~Renaud.
\newblock Quelques propri\'et\'es des applications analytiques d'une boule de
  dimension infinie dans une autre.
\newblock {\em Bull. Sci. Math. (2)}, 97:129--159 (1974), 1973.

\end{thebibliography}
\end{document}